\pdfoutput=1
\documentclass[10pt]{amsart}
\usepackage{amscd}
\usepackage{amssymb}
\usepackage[all]{xy}
\usepackage{graphicx}
\usepackage{hyperref}
\hypersetup{colorlinks=false}
\title[Rigidity, Residues and Duality]{Rigidity, Residues and Duality: Overview 
and Recent Progress}
\date{30 January 2021}
\usepackage{libertine} 
\usepackage[scaled=0.92]{roboto}  %{cabin} 
\usepackage[T1]{fontenc} 
\usepackage[scaled=0.87]{inconsolata}
\usepackage[libertine,varbb]{newtxmath}
\usepackage[cal=cm, calscaled=1.00, frak=euler, frakscaled=0.98]{mathalfa}
\author{Amnon Yekutieli}
\address{Department of  Mathematics,
Ben Gurion University, Be'er Sheva 84105, Israel}
\email{amyekut@math.bgu.ac.il}
\thanks{{\em Mathematics Subject Classification} 2020.
Primary: 14F08. Secondary: 18G80, 13D09, 14A20, 18F20, 16E45.}
\keywords{Grothendieck Duality, rigid dualizing complexes, 
rigid residue complexes,  schemes, algebraic spaces, 
Deligne-Mumford stacks, derived categories, derived functors, DG rings, 
DG modules.}
\thanks{Supported by the Israel Science Foundation grant no.\ 824/12.}

\newtheorem{thm}[equation]{Theorem}
\newtheorem{cor}[equation]{Corollary}

\newtheorem{thmstar}[equation]{Theorem$\boldsymbol{^{(*)}}$}
\newtheorem{corstar}[equation]{Corollary$\boldsymbol{^{(*)}}$}
\theoremstyle{definition}
\newtheorem{dfn}[equation]{Definition}
\newtheorem{rem}[equation]{Remark}
\newtheorem{exa}[equation]{Example}

\numberwithin{equation}{section}
\newcommand{\iso}{\xrightarrow{% 
\smash{\raisebox{-0.5ex}{\ensuremath{\scriptstyle \simeq  \mspace{2mu}}}}}} 
\renewcommand{\equiv}{\xrightarrow{% 
\smash{\raisebox{-0.5ex}{\ensuremath{\scriptstyle \approx  \mspace{2mu}}}}}} 
\newcommand{\inj}{\hookrightarrow}

\newcommand{\xar}{\xrightarrow}
\newcommand{\opn}{\operatorname}
\newcommand{\cat}[1]{\operatorname{\mathsf{#1}}}

\newcommand{\cd}{\,{\cdot}\,}

\newcommand{\ol}{\overline}
\newcommand{\rmitem}[1]{\item[\text{\textup{(#1)}}]}
\newcommand{\mfrak}[1]{\mathfrak{#1}}
\newcommand{\mcal}[1]{\mathcal{#1}}

\newcommand{\mbf}[1]{\mathbf{#1}}
\newcommand{\mrm}[1]{\mathrm{#1}}
\newcommand{\mbb}[1]{\mathbb{#1}}
\newcommand{\OO}{\mcal{O}}
\newcommand{\MM}{\mcal{M}}
\newcommand{\NN}{\mcal{N}}
\newcommand{\KK}{\mcal{K}}

\newcommand{\JJ}{\mcal{J}}
\newcommand{\II}{\mcal{I}}

\newcommand{\Ga}{\Gamma}
\newcommand{\si}{\sigma}
\newcommand{\la}{\lambda}

\renewcommand{\th}{\theta}
\newcommand{\om}{\omega}

\newcommand{\ga}{\gamma}

\newcommand{\Om}{\Omega}
\newcommand{\p}{\mfrak{p}}

\newcommand{\kk}{\bsym{k}}

\newcommand{\K}{\mathbb{K}}

\newcommand{\Z}{\mathbb{Z}}

\newcommand{\DR}{\mrm{R}}

\newcommand{\Hom}{\mcal{H}om}
\newcommand{\XX}{\mfrak{X}}
\newcommand{\YY}{\mfrak{Y}}

\newcommand{\tup}[1]{\textup{#1}}
\newcommand{\hsp}[1]{\hspace{#1 ex}}
\newcommand{\bsym}[1]{\boldsymbol{#1}}

\newcommand{\ot}{\otimes}

\newcommand{\til}[1]{\tilde{#1}}

\newcommand{\what}[1]{\widehat{#1}}

\newcommand{\sub}{\subseteq}

\newcommand{\lb}{\linebreak}

\newcommand{\twoiso}{%
\stackrel{\smash{\raisebox{-0.5ex}{\ensuremath{\scriptstyle \simeq%
\mspace{2mu}}}}}{\Longrightarrow}}
\newcommand{\eftover}{\mspace{0.0mu} /_{\mspace{-2mu} \mrm{eft}}\mspace{1.5mu}}
\newcommand{\ftover}{\mspace{0.0mu} /_{\mspace{-2mu} \mrm{ft}}\mspace{1.5mu}}
\renewcommand{\over}{\mspace{0mu} / \mspace{0mu}}

\newcommand{\scalesubsub}[1]{\scalebox{0.4}{$#1$}}

\begin{document}

\begin{abstract}
In this article we explain the theory of {\em rigid residue complexes} 
in commutative algebra and algebraic geometry, summarizing the background, 
recent results and anticipated future results. 
Unlike all previous approaches to Grothendieck Duality, 
the rigid approach concentrates on the construction of rigid residue complexes 
over rings, and their intricate yet robust properties. The geometrization, 
i.e.\ the passage to rigid residue complexes on schemes and 
Deligne-Mumford (DM) stacks, by gluing, is fairly easy. In the geometric part 
of the theory, the main results are the {\em Rigid Residue Theorem} and the 
{\em Rigid Duality Theorem} for proper maps between schemes, and for tame proper 
maps between DM stacks. 
\end{abstract}

\maketitle

\tableofcontents		

%\cleardoublepage
\setcounter{section}{-1}
\section{Introduction}

This article is about the {\em rigid approach to Grothendieck Duality} (GD).
We survey the main features of this approach, discussing both results that were 
already obtained, and future expected results. 
 
Throughout this article we work over a fixed base ring $\K$, which is a finite 
dimensional regular noetherian commutative ring (e.g.\ a field or the ring of 
integers $\Z$). All rings are commutative $\K$-rings, and by default they are 
essentially finite type (EFT) over $\K$.

In Section 1 we review the {\em squaring operation}. Given a ring $A$, the 
squaring operation is a functor $\opn{Sq}_{A / \K}$ from the derived category 
$\cat{D}(A)$ to itself. This is a quadratic functor, and -- unless $A$ is flat 
over $\K$ -- its construction requires the use of {\em commutative and 
noncommutative DG rings}; see \cite{Ye6}. The squaring operation has 
{\em backward functoriality} over any ring homomorphism $u : A \to B$, and 
{\em forward functoriality} when $u$ is {\em essentially \'etale}.

A {\em rigid complex} over $A$ (relative to $\K$) is a pair $(M, \rho)$, where 
$M \in \cat{D}^{\mrm{b}}_{\mrm{f}}(A)$, and 
$\rho : M \iso \opn{Sq}_{A / \K}(M)$ is an isomorphism in $\cat{D}(A)$. 
There is a pretty obvious notion of {\em rigid morphism} between rigid 
complexes (see Definition \ref{dfn:235} ). This notion extends to {\em rigid 
backward morphisms} over any ring homomorphism $u : A \to B$, and to {\em rigid 
forward morphisms} when $u$ is essentially \'etale.

A {\em rigid dualizing complex} over $A$ is a rigid complex 
$(R_A, \rho_A)$ as above, such that $R_A$ is a dualizing complex over $A$ in 
the sense of \cite{RD}. This concept is discussed in Section 2. 
The ring $A$ has a rigid dualizing complex
$(R_A, \rho_A)$, and it is unique up to a unique rigid isomorphism. 
If $u : A \to B$ is a finite ring homomorphism, then there is a unique 
nondegenerate rigid backward morphism 
$\opn{tr}^{\mrm{rig}}_u : (R_B, \rho_B) \to (R_A, \rho_A)$
called the {\em rigid trace morphism}. 
If $v : A \to A'$ is an essentially \'etale ring homomorphism, then 
there is a unique nondegenerate rigid forward morphism 
$\opn{q}^{\mrm{rig}}_v : (R_A, \rho_A) \to (R_{A'}, \rho_{A'})$
called the {\em rigid \'etale-localization morphism}. The morphisms 
$\opn{tr}^{\mrm{rig}}_u$ and $\opn{q}^{\mrm{rig}}_v$
commute with each other in the appropriate sense. 

Having the rigid dualizing complex $(R_A, \rho_A)$ on every EFT $\K$-ring $A$
allows to construct the {\em twisted induction pseudofunctor}, in a
totally algebraic way (rings only, no geometry), equipped with forward and 
backward functorialities. All the material up to this point (the content of 
Sections 1-2) is in the papers \cite{Ye6} and \cite{OSY}, the latter joint 
with M. Ornaghi and S. Singh. 

Looking to the future, we plan to study a more refined object than the rigid 
dualizing complex: it is the {\em rigid residue complex} 
$(\KK_A, \rho_A)$ of a ring $A$. Given an essentially \'etale ring homomorphism 
$v : A \to A'$, there is the {\em nondegenerate rigid \'etale-localization 
homomorphism} $\opn{q}^{\mrm{res}}_{v} : \KK_{A} \to \KK_{A'}$. 
For a finite ring homomorphism $u : A \to B$ there is the {\em nondegenerate 
rigid trace homomorphism} 
$\opn{tr}^{\mrm{res}}_{u} : \KK_{B} \to \KK_{A}$. 
The latter extends to a functorial backward homomorphism 
$\opn{tr}^{\mrm{ires}}_{u} : \KK_B \to \KK_A$ 
for an arbitrary ring homomorphism $u$, called the {\em ind-rigid trace}, but 
$\opn{tr}^{\mrm{ires}}_{u}$ is only a homomorphism of graded $A$-modules.
We should stress that $\opn{q}^{\mrm{res}}_{v}$ and 
$\opn{tr}^{\mrm{res}}_{v}$ are actual homomorphisms of complexes (see Remark 
\ref{rem:295} regarding notation). 
This topic, summarized in Section 3, is planned for the paper \cite{Ye8}. 

Rigid residue complexes can be easily glued on finite type (FT) $\K$-schemes, 
and they still have the ind-rigid trace 
$\opn{tr}^{\mrm{ires}}_f : f_{*}(\KK_Y) \to \KK_X$
for an arbitrary scheme map $f : Y \to X$. The twisted induction pseudofunctor 
from Section 2 now becomes the geometric {\em  twisted inverse image 
pseudofunctor} $f \mapsto f^!$. We expect to prove the Rigid Residue Theorem
\ref{thm:1} and the Rigid Duality Theorem \ref{thm:2} for proper maps of EFT 
$\K$-schemes, thus recovering almost all the results from the original book 
“Residues and Duality” \cite{RD}, yet in a very explicit way.
This topic is covered in Section 4, and planned for the paper \cite{Ye9}. 

For every scheme $X$ the assignment $X' \mapsto \KK_{X'}$ 
is a quasi-coherent sheaf on the small \'etale site $X_{\rm{et}}$.
This implies that every finite type Deligne-Mumford (DM) $\K$-stack $\XX$ 
admits a rigid residue complex $\KK_{\XX}$. Here too we have the 
twisted inverse image pseudofunctor $f \mapsto f^!$ . For a map 
$f : \YY \to \XX$ of FT DM $\K$-stacks there is the ind-rigid trace
$\opn{tr}^{\mrm{ires}}_f : f_{*}(\KK_Y) \to \KK_X$. Under a mild 
technical condition, we expect to prove the Rigid Residue Theorem \ref{thm:10}
for proper maps of DM stacks, and the Rigid Duality Theorem \ref{thm:11} for 
such maps that are also tame. This is the content of Section 5 of the article, 
and will be written in the paper \cite{Ye10}. 

To finish the Introduction, a few words on the history of the rigid approach to 
GD. Rigid dualizing complexes were invented by M. Van den Bergh in the 1990's, 
in the noncommutative setting. This concept was imported to commutative 
algebra by J.J. Zhang and the author around 2005, where it was made functorial, 
and it was also expanded to the arithmetic setting (namely allowing a base 
ring $\K$ that is not a field). More details on the history of the topic, with 
references, can be found  in Remark \ref{rem:279}. 
Many of the ideas in this article were already mentioned in our lectures from 
2013, see \cite{Ye5}.

\medskip \noindent 
{\em Acknowledgments.}
Because this is a project spanning almost 15 years (and counting), it is 
impossible to thank all the people who had contributed to it. Therefore I shall 
only mention my collaborator James J. Zhang with whom the project began; my 
current collaborators Mattia Ornaghi and Saurabh Singh;   
Martin Olsson and Dan Edidin who assisted me in my attempts to understand 
DM stacks; and Jack Hall and Amnon Neeman who explained some of their work to 
me.  The project is currently supported by the Israel Science Foundation 
grant no.\ 824/18.

% \cleardoublepage
\section{The Squaring Operation and Rigid Complexes over Rings}
 
We fix a base ring $\K$, which is regular noetherian of finite Krull 
dimension (e.g.\ a field or $\Z$).
Recall that a commutative $\K$-ring $A$ is called {\em essentially finite type} 
(EFT) if it is a localization of a finite type $\K$-ring. 
Note that such a ring $A$ is noetherian and of finite Krull dimension.
The category of EFT commutative $\K$-rings is denoted by 
$\cat{Rng} \eftover \K$.
All rings in this paper belong to $\cat{Rng} \eftover \K$.

For a ring $A$ we denote by $\cat{M}(A)$ the abelian category of 
$A$-modules, and by $\cat{C}(A)$ the DG category of complexes of 
$A$-modules. The strict subcategory 
$\cat{C}_{\mrm{str}}(A) \sub \cat{C}(A)$
has the same objects, but its morphisms are the degree $0$ cocycles (i.e.\ 
the degree $0$ homomorphisms $\phi : M \to  N$ that commute with the 
differentials). The derived category is $\cat{D}(A)$.
The categorical localization functor 
$\opn{Q} : \cat{C}_{\mrm{str}}(A) \to \cat{D}(A)$
is the identity on objects, and it is universal for inverting 
quasi-isomorphisms in $\cat{C}_{\mrm{str}}(A)$. 
Note that $\cat{C}_{\mrm{str}}(A)$ is an abelian category, $\cat{D}(A)$ is a 
triangulated category, and the functor $\opn{Q}$ is $A$-linear. 
Inside $\cat{D}(A)$ there is the full triangulated subcategory
$\cat{D}^{\mrm{b}}_{\mrm{f}}(A)$
of complexes with bounded finitely generated cohomology. 
All the necessary facts about derived categories can be found in our book 
\cite{Ye7}. Furthermore, the contents of Sections 1-3 are treated in detail in 
this book, under the assumption that the rings in question are flat over $\K$. 

Following \cite{Ye6} and \cite[Chapter 3]{Ye7}, by {\em commutative DG 
$\K$-ring} we mean a nonpositive strongly commutative DG ring 
$\til{A} = \bigoplus_{i \leq 0} \til{A}^i$,
equipped with a ring homomorphism $\K \to \til{A}^0$.
Strong commutativity means that $b \cd a = (-1)^{i \cd j} \cd a \cd b$ for 
all $a \in \til{A}^i$ and $b \in \til{A}^j$, and that $a \cd a = 0$ when $i$ is 
odd. The DG ring $\til{A}$ is said to be K-flat if it is K-flat as a DG 
$\K$-module. 

Let $A$ be a ring. A {\em K-flat commutative DG ring resolution} of $A$ is  
a quasi-isomorphism $\til{A} \to A$ from a K-flat commutative DG $\K$-ring 
$\til{A}$. Such resolutions always exist. 

Here is the main theorem of \cite{Ye6}:

\begin{thm}[\cite{YZ1}, \cite{Ye6}] \label{thm:230}
There is a functor
\[ \opn{Sq}_{A / \K} : \cat{D}(A) \to \cat{D}(A) \]
called the {\em squaring operation}, equipped with a canonical isomorphism 
\[ \opn{Sq}_{A / \K}(M) \cong  
\opn{RHom}_{\til{A} \ot_{\K} \til{A}}(A, M \otimes^{\mrm{L}}_{\K} M) \]
for every K-flat  commutative DG ring resolution $\til{A} \to A$. 
\end{thm}

See \cite[Theorem 0.3.4]{Ye6} for the precise statement. It is interesting to 
mention that the proof requires the use of {\em noncommutative DG rings}. 

\begin{thm}[\cite{YZ1}, \cite{Ye6}] \label{thm:231}
The functor $\opn{Sq}_{A / \K}$ is a {\em quadratic functor}, namely for a 
morphism $\phi : M \to N$ in $\cat{D}(A)$ and an element $a \in A$, there is 
equality 
\[ \opn{Sq}_{A / \K} (a \cd \phi) = a^2 \cd \opn{Sq}_{A / \K} (\phi) . \]
\end{thm}

\begin{dfn}[\cite{YZ1}, \cite{Ye6}] \label{dfn:230}
A {\em rigid complex} over $A$ relative to $\K$ is a pair 
$(M, \rho)$, consisting of a complex 
$M \in \cat{D}^{\mrm{b}}_{\mrm{f}}(A)$, and an isomorphism
$\rho : M \iso \opn{Sq}_{A / \K}(M)$ in $\cat{D}(A)$,
called a {\em rigidifying isomorphism} for $M$. 
\end{dfn}

\begin{dfn}[\cite{YZ1}, \cite{Ye6}] \label{dfn:231}
Suppose $(M, \rho)$ and $(N, \sigma)$ are rigid complexes. A {\em 
rigid morphism} $\phi : (N, \sigma) \to (M, \rho)$
is a morphism $\phi : N \to M$ in $\cat{D}(A)$, such that the diagram
\[ \UseTips \xymatrix @C=7ex @R=4ex {
N
\ar[r]^(0.34){\sigma}
\ar[d]_{\phi}
&
\opn{Sq}_{A / \K}(N)
\ar[d]^{\opn{Sq}_{A / \K}(\phi)}
\\
M
\ar[r]^(0.34){\rho}
&
\opn{Sq}_{A / \K}(M)
} \]
in $\cat{D}(A)$ is commutative.

We denote by $\cat{D}(A)_{\mrm{rig} / \K}$
the category of rigid complexes and rigid morphisms between them.
\end{dfn}

\begin{dfn}[\cite{Ye7}] \label{dfn:232}
A complex $M \in \cat{D}^{\mrm{b}}_{\mrm{f}}(A)$
is said to have the {\em derived Morita property} if the derived homothety 
morphism $A \to \opn{RHom}_{A}(M, M)$ in $\cat{D}(A)$ is an isomorphism.
\end{dfn}

Dualizing complexes (to be recalled later) and tilting complexes (see 
\cite[Chapter 14]{Ye7}) have the derived Morita property. 

\begin{thm}[\cite{YZ1}, \cite{Ye6}] \label{thm:233}
If $(M, \rho)$ is a rigid complex, such that $M$ has the derived Morita 
property, then the only  automorphism of $(M, \rho)$ in
$\cat{D}(A)_{\mrm{rig} / \K}$ is the identity. 
\end{thm}

The key to the proof of this theorem is the quadratic property of the functor 
$\opn{Sq}_{A / \K}$, see Theorem \ref{thm:231}.

The squaring operation $\opn{Sq}_{A / \K}$  has two kinds of functorialities, 
which we describe next.  First there is {\em backward functoriality}. 
Given a ring homomorphism $u : A \to B$, and
complexes $M \in \cat{D}(A)$ and $N \in \cat{D}(B)$, a morphism 
$\th : N \to  M$ in $\cat{D}(A)$ is called a {\em backward morphism} over $u$. 
The backward morphism $\th$ is called {\em nondegenerate} if the morphism
$N \to \opn{RHom}_A(B, M)$ in $\cat{D}(B)$,
which corresponds to $\th$ by adjunction, is an isomorphism.
 
The next theorem extends Theorem \ref{thm:230}. 

\begin{thm}[\cite{YZ1}, \cite{Ye6}] \label{thm:234}
Let $u : A \to B$ be a ring homomorphism, and let 
$\th : N \to  M$ be a backward morphism in $\cat{D}(A)$ over $u$. Then 
there is a  backward morphism
\[ \opn{Sq}_{u / \K}(\th) :\opn{Sq}_{B / \K}(N) \to \opn{Sq}_{A / \K}(M) \]
in $\cat{D}(A)$ over $u$, with an explicit formula given suitable DG ring 
resolutions $\til{A} \to A$ and $\til{B} \to B$. 
The operation $\opn{Sq}_{u / \K}(\th)$ is functorial in $u$ and $\th$. 
\end{thm}

For a precise statement see \cite[Theorem 0.3.5]{Ye6}.

\begin{dfn}[\cite{YZ1}, \cite{OSY}] \label{dfn:233}
Let $u : A \to B$ be a ring homomorphism, let 
$(M, \rho) \in \cat{D}(A)_{\mrm{rig} / \K}$, and let 
$(N, \si) \in \cat{D}(B)_{\mrm{rig} / \K}$. 
A {\em rigid backward morphism} 
$\th :  (N ,\si) \to (M ,\rho)$
over $u$ is a backward morphism $\th : N \to M$ in $\cat{D}(A)$ over $u$,  
such that the diagram 
\[ \UseTips \xymatrix @C=7ex @R=4ex {
N
\ar[d]_{\th}
\ar[r]^(0.34){\si}
&
\opn{Sq}_{B / \K}(N)
\ar[d]^{\opn{Sq}_{u / \K}(\th)}
\\
M
\ar[r]^(0.34){\rho}
&
\opn{Sq}_{A / \K}(M)
} \]
in $\cat{D}(A)$ is commutative.
\end{dfn}

Next is a generalization of Theorem \ref{thm:233}.

\begin{thm}[\cite{YZ1}, \cite{OSY}] \label{thm:235}
In the situation of Definition \ref{dfn:233}, if $N \in \cat{D}(B)$ has the 
derived Morita property, then there is at most 
one nondegenerate rigid backward morphism 
$\th :  (N ,\si) \to (M ,\rho)$ over $u$. 
\end{thm}

The second functoriality of the squaring operation is {\em forward 
functoriality}. A homomorphism $v : A \to A'$ is called 
{\em essentially \'etale} if it is EFT and formally \'etale. 
We know that an essentially \'etale  homomorphism is flat. 
Given an essentially \'etale  homomorphism $v : A \to A'$, and complexes 
$M \in \cat{D}(A)$ and $M' \in \cat{D}(A')$, a morphism 
$\la : M \to  M'$ in $\cat{D}(A)$ is called a {\em forward morphism} over $v$. 
The forward morphism $\la$ is called {\em nondegenerate} if the morphism
$A' \ot_A M \to M'$ in $\cat{D}(A')$, which corresponds to $\la$ by 
adjunction, is an isomorphism.

\begin{thm}[\cite{YZ1}, \cite{OSY}] \label{thm:236}
Let $v : A \to A'$ be an essentially \'etale ring homomorphism, and let 
$\la : M \to  M'$ be a forward morphism in $\cat{D}(A)$ over $v$. Then 
there is a forward morphism 
\[ \opn{Sq}_{u / \K}(\la) :\opn{Sq}_{A / \K}(M) \to \opn{Sq}_{A' / \K}(M') \]
in $\cat{D}(A)$ over $v$, with an explicit formula given suitable DG ring 
resolutions $\til{A} \to A$ and $\til{B} \to B$. 
The operation $\opn{Sq}_{u / \K}(\la)$ is functorial in $v$ and $\la$. 
\end{thm}

In \cite{YZ1} we only considered the case of a {\em localization} homomorphism 
$v : A \to A'$. The essentially \'etale case, proved in \cite{OSY}, is much 
more difficult, and it relies on a detailed study of the diagonal embedding of 
$\opn{Spec}(A')$ in $\opn{Spec}(A' \ot_A A')$. 

\begin{dfn}[\cite{YZ1}, \cite{OSY}] \label{dfn:235} 
Let $v : A \to A'$ be an essentially \'etale ring homomorphism, let 
$(M, \rho) \in \cat{D}(A)_{\mrm{rig} / \K}$, and let 
$(M', \rho') \in \cat{D}(A')_{\mrm{rig} / \K}$. 
A {\em rigid forward morphism} 
$\la :  (M ,\rho) \to (M' ,\rho')$
over $v$ is a forward morphism $\la : M \to M'$ in $\cat{D}(A)$ over $v$,  
such that the diagram 
\[ \UseTips \xymatrix @C=7ex @R=4ex {
M
\ar[d]_{\la}
\ar[r]^(0.34){\rho}
&
\opn{Sq}_{B / \K}(M)
\ar[d]^{\opn{Sq}_{v / \K}(\la)}
\\
M'
\ar[r]^(0.34){\rho'}
&
\opn{Sq}_{A / \K}(M')
} \]
in $\cat{D}(A)$ is commutative.
\end{dfn}

Here is a second generalization of Theorem \ref{thm:233}.

\begin{thm}[\cite{OSY}] \label{thm:237}
In the situation of Definition \tup{\ref{dfn:235}}, if $M' \in \cat{D}(A')$ has 
the derived Morita property, then there is at most 
one nondegenerate rigid forward morphism 
$\la :  (M ,\rho) \to (M' ,\rho')$ over $v$. 
\end{thm}
 
The backward and forward squaring operations commute with each other. This is 
the content of the next theorem. 

\begin{thm}[\cite{OSY}] \label{thm:238}
Suppose we are given ring homomorphisms $u : A \to  B$ and $v : A \to  A'$, 
such 
that $v$ is essentially \'etale, and a backward morphism $\th : N \to M$ in 
$\cat{D}(A)$ over $u$. Define the ring $B' := A' \ot_A B$, and the complexes 
$M' := A' \ot_A M \in \cat{D}(A')$ and
$N' := B' \ot_B N \in \cat{D}(B')$. Let 
$\th' : N' \to M'$ be the induced backward morphism in $\cat{D}(A')$,
and let $\opn{q}_{M} : M \to M'$ and $\opn{q}_{N} : N \to N'$ be the canonical 
nondegenerate forward morphisms. We get the following commutative diagrams
\[ \UseTips \xymatrix @C=6ex @R=5ex {
A
\ar[r]^{u}
\ar[d]_{v}
&
B
\ar[d]^{w}
\\
A'
\ar[r]^{u'}
&
B'
}
\qquad \qquad 
\UseTips \xymatrix @C=6ex @R=5ex {
M
\ar[d]_{\opn{q}_{M}}
&
N
\ar[l]_{\th}
\ar[d]^{\opn{q}_{N}}
\\
M'
&
N'
\ar[l]_{\th'}
} \]
in $\cat{Rng} \eftover \K$ and in $\cat{D}(A)$. Then the diagram  
\begin{equation} \label{eqn:219}
\UseTips \xymatrix @C=10ex @R=6ex {
\opn{Sq}_{A / \K}(M)
\ar[d]_{\opn{Sq}_{v / \K}(\opn{q}_{M})}
&
\opn{Sq}_{B / \K}(N)
\ar[l]_{\opn{Sq}_{u / \K}(\th)}
\ar[d]^{\opn{Sq}_{w / \K}(\opn{q}_{N})}
\\
\opn{Sq}_{A / \K}(M')
&
\opn{Sq}_{B / \K}(N')
\ar[l]_{\opn{Sq}_{u' / \K}(\th')}
}
\end{equation}
in $\cat{D}(A)$ is commutative. 
\end{thm}

This theorem implies that there is no conflict between Definitions  
\ref{dfn:233} and \ref{dfn:235} when $A = B = A'$ and $u = v = \opn{id}_A$, in 
which case the morphisms $\th$ and $\la$ in these definitions are 
both forward and backward morphisms.

%\cleardoublepage
\section{Rigid Dualizing Complexes over Rings} \label{sec:RRD-rings}
 
Again $\K$ is a regular base ring, and $A$ is an EFT $\K$-ring. 
The next definition is taken from \cite{RD}, except that the derived Morita 
property (Definition \ref{dfn:232}) had not yet been introduced when 
\cite{RD} was written. 

\begin{dfn} \label{dfn:240}
A complex $R \in \cat{D}^{\mrm{b}}_{\mrm{f}}(A)$
is called {\em dualizing} if it has finite injective dimension and 
the derived Morita property.
\end{dfn}

\begin{dfn}[\cite{VdB}, \cite{YZ2}, \cite{Ye7}, \cite{OSY}] \label{dfn:241}
A {\em rigid dualizing complex} over $A$ relative to $\K$ is a rigid 
complex $(R, \rho)$, as in Definition \ref{dfn:230}, such that $R$ is a 
dualizing complex. 
\end{dfn}

\begin{thm}[{\cite{VdB}, \cite{YZ2}, \cite{OSY}}] \label{thm:120}
Let $A$ be an EFT $\K$-ring. The ring $A$ has a rigid dualizing complex 
$(R_A, \rho_A)$, and it is unique, up to a unique isomorphism in 
$\cat{D}(A)_{\mrm{rig} / \K}$. 
\end{thm}

The uniqueness part of this theorem is a combination of a result of 
Grothendieck (see \cite[Theorem V.3.1]{RD} or \cite[Theorem 13.1.35]{Ye7}), 
with a variation of Theorem \ref{thm:233} (see \cite[Theorem 13.5.4]{Ye7}).
Existence is much harder. The strategy of proof in \cite{OSY} is this: we
factor the structure homomorphism $u : \K \to A$ into 
$u = u_{\mrm{loc}} \circ u_{\mrm{fin}} \circ u_{\mrm{pl}}$,
where $u_{\mrm{pl}} : \K \to A_{\mrm{pl}}$ is a homomorphism to a polynomial 
ring over $\K$; 
$u_{\mrm{fin}} : A_{\mrm{pl}} \to A_{\mrm{ft}}$ is a surjection;
and $u_{\mrm{loc}} : A_{\mrm{ft}} \to A$ is a localization homomorphism. 
Then we show that the dualizing complexes $R_{\mrm{pl}}$,  $R_{\mrm{ft}}$ and 
$R$, over the rings $A_{\mrm{pl}}$,  $A_{\mrm{ft}}$ and $A$ respectively, that 
are constructed in the course of the proof of \cite[Theorem 13.1.34]{Ye7}, are 
all rigid. This is achieved using {\em induced rigidity} for essentially smooth 
ring homomorphisms, and {\em coinduced rigidity} for finite ring homomorphisms. 

\begin{thm}[\cite{YZ2}, \cite{OSY}] \label{thm:241} 
Let $u : A \to B$ be a finite ring homomorphism.
There exists a unique nondegenerate rigid backward morphism
\[ \opn{tr}^{\mrm{rig}}_u : (R_B, \rho_B) \to (R_A, \rho_A) \]
in $\cat{D}(A)$ over $u$, called the {\em rigid trace morphism}.
\end{thm}

\begin{cor}[\cite{YZ2}, \cite{OSY}] \label{cor:241}
Let $A \xar{u} B \xar{v} C$ be finite ring homomorphisms. Then 
$\opn{tr}^{\mrm{rig}}_u \circ \opn{tr}^{\mrm{rig}}_v = 
\opn{tr}^{\mrm{rig}}_{v \circ u}$, as backward morphisms 
$R_C \to R_A$ in $\cat{D}(A)$ over $v \circ u$.
\end{cor}

\begin{thm}[\cite{YZ2}, \cite{OSY}] \label{thm:242} 
Let $v : A \to A'$ be an essentially \'etale ring homomorphism.
There exists a unique nondegenerate rigid forward morphism
\[ \opn{q}^{\mrm{rig}}_v : (R_A, \rho_A) \to (R_{A'}, \rho_{A'}) \]
in $\cat{D}(A)$ over $v$, called the {\em rigid \'etale-localization morphism}.
\end{thm}

\begin{cor}[\cite{YZ2}, \cite{OSY}] \label{cor:242}
Let $A \xar{v} A' \xar{v'} A''$ be essentially \'etale ring homomorphisms.
Then $\opn{q}^{\mrm{rig}}_{v'} \circ \opn{q}^{\mrm{rig}}_v = 
\opn{q}^{\mrm{rig}}_{v' \circ v}$, as forward morphisms 
$R_A \to R_{A''}$ in $\cat{D}(A)$ over $v' \circ v$.
\end{cor}

\begin{thm}[\cite{YZ2}, \cite{OSY}] \label{thm:240} 
Suppose $u : A \to  B$ and $v : A \to  A'$ are ring homomorphisms, with $u$
finite and $v$ essentially \'etale. Define the ring $B' := A' \ot_A B$.
We get the following commutative diagram
\[ \UseTips \xymatrix @C=6ex @R=5ex {
A
\ar[r]^{u}
\ar[d]_{v}
&
B
\ar[d]^{w}
\\
A'
\ar[r]^{u'}
&
B'
} \]
in $\cat{Rng} \eftover \K$, in which $u'$ is finite and $w$ is essentially 
\'etale. Then the diagram  
\[ \UseTips \xymatrix @C=10ex @R=6ex {
R_A
\ar[d]_{\opn{q}^{\mrm{rig}}_{v}}
&
R_B
\ar[l]_{\opn{tr}^{\mrm{rig}}_u}
\ar[d]^{\opn{q}^{\mrm{rig}}_{w}}
\\
R_{A'}
&
R_{B'}
\ar[l]_{\opn{tr}^{\mrm{rig}}_{u'}}
} \]
in $\cat{D}(A)$ is commutative. 
\end{thm}

For a ring $A$ with rigid dualizing complex $(R_A, \rho_A)$, let us define the 
{\em rigid auto-duality functor}
\begin{equation} \label{eqn:130}
\opn{D}^{\mrm{rig}}_A := \opn{RHom}_{A}(-, R_A) : \cat{D}(A) \to \cat{D}(A) . 
\end{equation}
According to \cite{RD}, or \cite[Chapter 13]{Ye7}, the functor 
$\opn{D}^{\mrm{rig}}_A$ induces an equivalence
$\opn{D}^{\mrm{rig}}_A : \cat{D}^{\star}_{\mrm{f}}(A)^{\mrm{op}} \to 
\cat{D}^{-\star}_{\mrm{f}}(A)$,
where $\star$ and $-\star$ are reverse boundedness conditions. 
In particular there is the {\em evaluation isomorphism} 
\begin{equation} \label{eqn:245}
\opn{ev} : \opn{Id} \iso \opn{D}^{\mrm{rig}}_A \circ \opn{D}^{\mrm{rig}}_A 
\end{equation}
of triangulated functors from $\cat{D}^{+}_{\mrm{f}}(A)$ to itself.

The next discussion involves {\em $2$-categories} and {\em pseudofunctors}.
The only textbook reference for $2$-categories seems to be in 
\cite[Chapter XII]{Mc}. There is a brief discussion of $2$-categories and 
pseudofunctors in 
\cite[Section \href{https://stacks.math.columbia.edu/tag/003G}
{\tt tag=003G}]{SP}, 
and a summary of $2$-categorical notation in Section 8.1 of 
our book \cite{Ye7}. A detailed exposition of $2$-categories and pseudofunctors
can be found in Sections 1-2 of our paper \cite{Ye4}. We should stress that our 
$2$-categories are strict, and our pseudofunctors are normalized. 
 
Recall that to a ring homomorphism $u : A \to B$ we can associate the {\em 
derived induction functor}
\begin{equation} \label{eqn:220}
\mrm{LInd}_u  := B \ot^{\mrm{L}}_{A} (-) : \cat{D}^{-}_{\mrm{f}}(A) \to
\cat{D}^{-}_{\mrm{f}}(B)  . 
\end{equation}
It is the left derived functor of the {\em induction functor} 
$\mrm{Ind}_u  := B \ot_{A} (-) : \cat{M}(A) \to \cat{M}(B)$.

Consider the $2$-category of $\K$-linear triangulated categories, which we 
denote by \lb $\cat{TrCat} \over \K$.
Derived induction is actually a pseudofunctor
\begin{equation} \label{eqn:221}
\mrm{LInd} :  \cat{Rng} \eftover \K \to \cat{TrCat} \over \K ,
\end{equation}
sending a ring $A$ to the category $\cat{D}^{-}_{\mrm{f}}(A)$, and 
a ring homomorphism $u : A \to B$ to the functor 
$\mrm{LInd}_u$. The composition isomorphisms of $\mrm{LInd}$ arise from those 
of the plain induction pseudofunctor $\opn{Ind}$; and 
these are just the associativity isomorphisms for tensor products. 

The next result is our ring-theoretic version of the geometric {\em 
twisted inverse pseudofunctor} $f \mapsto f^!$ from \cite{RD}. 
See Theorem \ref{thm:290}  for our geometric result. 

\begin{thm}[{\cite{YZ2}, \cite{OSY}}] \label{thm:121}
There is a unique pseudofunctor 
\[ \opn{TwInd} : \cat{Rng} \eftover \K \to \cat{TrCat} \over \K \]
called {\em twisted induction}, with these properties\hsp{0.1}\tup{:}
\begin{enumerate} 
\item To an object $A \in \cat{Rng} \eftover \K$ it assigns the triangulated 
category $\cat{D}^{+}_{\mrm{f}}(A)$.
 
\item To a morphism $u : A \to B$ in $\cat{Rng} \eftover \K$ it assigns the 
triangulated functor 
\[ \opn{TwInd}_{u} \ := \ 
\opn{D}_{B}^{\mrm{rig}} \circ \opn{LInd}_u \circ \opn{D}_{A}^{\mrm{rig}}
\ :  \ \cat{D}^{+}_{\mrm{f}}(A) \to \cat{D}^{+}_{\mrm{f}}(B) . \]
 
\item The composition isomorphisms of \hsp{0.1} $\opn{TwInd}$ come from the 
evaluation isomorphisms of the rigid auto-duality functors 
\tup{(\ref{eqn:245})}, combined with the composition isomorphisms of the 
pseudofunctor $\opn{LInd}$.
\end{enumerate}
\end{thm}

In the next section we shall see how $\opn{TwInd}$ interacts with rigid traces 
and rigid \'etale-localizations.

\begin{rem} \label{rem:279}
Observe that the pseudofunctor $\opn{TwInd}$ is constructed {\em purely by 
algebraic methods}, namely without geometry, and especially without {\em 
global duality}, which is a key ingredient in all earlier approaches to GD. 
The only other such purely algebraic construction we are aware of is the very 
recent one by P. Scholze and D. Clausen \cite{Sc}. They use {\em condensed 
mathematics}, and so far it works only for FT $\Z$-rings. 
\end{rem}

\begin{rem} \label{rem:278}
Here are a few historical details on the squaring operation and 
rigid dualizing complexes. 

{\em Noncommutative dualizing complexes} were introduced in our 1992 paper
\cite{Ye3}. These are complexes of bimodules over a NC ring $A$, central over 
a base field $\K$, satisfying conditions similar to those of commutative 
dualizing complexes. 
 
{\em Noncommutative rigid dualizing complexes} were invented by M. Van 
den Bergh in 1997, in his seminal paper \cite{VdB}. Observe the similarity 
between the squaring operation (Theorem \ref{thm:230}) and {\em Hochschild 
cohomology}. This is not a coincidence. An idea similar to that of a NC rigid 
dualizing complex was independently discovered by M. Kontsevich in the late 
1990's, and it is one of the cornerstones of the theory of {\em Calabi-Yau 
categories} that emerged from the work of Kontsevich and others. These NC 
theories are explained in \cite[Chapter 18]{Ye7}.
 
J.J. Zhang and the author imported the concept of rigid dualizing complexes 
back to commutative algebra and algebraic geometry. We added two features 
to the original concept of NC rigidity of Van den Bergh: the {\em functoriality 
of rigid complexes}, and the {\em passage from base field to base ring}, which 
is sometimes referred to as the {\em arithmetic setting}. This 
was done in our papers \cite{YZ1} and \cite{YZ2} from around 2008. 
As already mentioned in Section 1, the generalization from base field to base 
ring required the use of {\em DG ring resolutions}. 

Unfortunately there were several serious errors regarding the manipulation of 
DG rings in the paper \cite{YZ1}, and these errors affected some constructions 
and proofs in \cite{YZ1} and \cite{YZ2}. Interestingly, our main results in 
these two paper were correct, as we now know. 

The errors in \cite{YZ1} and \cite{YZ2} were discovered by the authors of 
\cite{AILN} in 2010, and they also fixed one of these errors. 
(However, contrary to the norms of ethical scientific conduct, the authors of 
\cite{AILN} failed to inform us of the errors in our work before publishing 
their own paper.) 
 
This accident --  the errors regarding DG rings -- took a lot of time and 
effort to correct. Some of the errors were fixed in our paper \cite{Ye6} from  
2016 , and the rest are corrected in the upcoming paper \cite{OSY}, which is 
joint work with M. Ornaghi and S. Singh.  
This is one of the main reasons the project discussed in this article was 
delayed for so many years.
\end{rem}

%\cleardoublepage
\section{Rigid Residue Complexes over Rings} \label{sec:rrr-rings}

Now we enter the zone of "anticipated results": many of the results still
only have sketchy proofs, and should therefore be regarded with caution. Such a 
theorem will be labeled {\bf Theorem$^{\bsym{(*)}}$}.

Recall that we are still working in $\cat{Rng} \eftover \K$, where $\K$ is a 
regular noetherian base ring.
If $L \in \cat{Rng} \eftover \K$ is a field, then its rigid dualizing 
complex $R_L$ must be isomorphic to $L[d]$ for some integer $d$. We define the 
{\em rigid dimension} of $L$ to be
$\opn{rig{.}dim}_{\K}(L) := d$.

\begin{exa}
If the base ring $\K$ is a field, then 
$\opn{rig{.}dim}_{\K}(L) = \opn{tr{.}deg}_{\K}(L)$, the transcendence degree.
On the other hand, for $\K = \Z$ and $L = \mbb{F}_p$ we have
$\opn{rig{.}dim}_{\Z}(\mbb{F}_p) = -1$.
\end{exa}

Now take an arbitrary $A \in \cat{Rng} \eftover \K$. 
For a prime ideal $\p \sub A$, with residue field $\bsym{k}(\p)$,
we define 
$\opn{rig{.}dim}_{\K}(\p) := \opn{rig{.}dim}_{\K}(\bsym{k}(\p))$.
The resulting function
$\opn{rig{.}dim}_{\K} : \opn{Spec}(A) \to \Z$
has the expected property: it drops by $1$ if $\p \sub \mfrak{q}$ is an 
immediate specialization of primes. 
The function $\opn{rig{.}dim}_{\K}$ is bounded, since 
$\opn{Spec}(A)$ is finite dimensional. 
See \cite[Theorem 13.3.3]{Ye7}.

For any $\p \in \opn{Spec}(A)$ we denote by $J(\p)$ the injective hull of the
$A_{\p}$-module $\bsym{k}(\p)$. This is an indecomposable injective $A$-module. 
According to the Matlis classification, every injective $A$-module is a direct 
sum, often infinite, of such modules $J(\p)$.
See \cite[Section 13.2]{Ye7} for a detailed discussion, and Remark 
\ref{rem:260} on the possible functorial properties of $J(\p)$.  

The following  definition is a refinement of the 
concept of {\em residual complex} from \cite{RD}. 

\begin{dfn} \label{dfn:250}
A {\em rigid residue complex} over $A$ relative to $\K$ is a rigid
dualizing complex {\em $(\mcal{K}_A, \rho_A)$} of this special sort: 
for every $d \in \Z$ there is an isomorphism of $A$-modules
\[ \mcal{K}_A^{-d} \cong 
\bigoplus_{\substack{ \p \, \in \, \opn{Spec}(A) \\  
\opn{rig{.}dim}_{\K} \mspace{-2mu}(\p) \, = \, d }}  J(\p) \ . \]
\end{dfn}
 
Recall the category of rigid complexes $\cat{D}(A)_{\mrm{rig} / \K}$ from 
Definition \ref{dfn:231}.

\begin{dfn} \label{dfn:251}
A {\em morphism of rigid residue complexes} 
$\phi : (\mcal{K}_A, \rho_A) \to (\mcal{K}'_A, \rho'_A)$
is a homomorphism 
$\phi : \mcal{K}_A \to \mcal{K}'_A$ in $\cat{C}_{\mrm{str}}(A)$, 
such that 
$\opn{Q}(\phi) : (\mcal{K}_A, \rho_A) \to (\mcal{K}'_A, \rho'_A)$
is a morphism in $\cat{D}(A)_{\mrm{rig} / \K}$.
We denote by $\cat{C}(A)_{\mrm{res} / \K}$
the category of rigid residue complexes. 
\end{dfn}

Note the unusual {\em mixed nature} of this definition: a morphism $\phi$ 
in $\cat{C}(A)_{\mrm{res} / \K}$ is literally is homomorphism of complexes, but 
the condition it must satisfy is in the derived category.  

\begin{thm}[\cite{Ye7}, \cite{OSY}] \label{thm:125}
Let $A$ be an EFT $\K$-ring. The ring $A$ has a rigid residue complex 
$(\KK_A, \rho_A)$ and it is unique, up to a unique isomorphism in 
$\cat{C}(A)_{\mrm{res} / \K}$. 
\end{thm}

This theorem is not actually stated in either of these references, yet it is 
easily proved by combining several results from them (Theorem \ref{thm:120} 
above, proved in \cite{OSY}, and \cite[Theorems 13.3.15 and 13.3.17]{Ye7}).

Let us now mention several important functorial properties of 
rigid residue complexes.

\begin{thmstar}[\cite{Ye8}] \label{thm:126}
 Suppose $v : A \to A'$ is an essentially \'etale homomorphism of EFT 
$\K$-rings. There is a unique nondegenerate forward homomorphism 
$\opn{q}^{\mrm{res}}_{v} : \mcal{K}_A \to \mcal{K}_{A'}$
in $\cat{C}_{\mrm{str}}(A)$, called the {\em rigid \'etale-localization 
homomorphism}, such that 
$\opn{Q}(\opn{q}^{\mrm{res}}_{v}) : 
(\KK_A, \rho_A) \to (\KK_{A'}, \rho_{A'})$
is the rigid \'etale-localization morphism $\opn{q}^{\mrm{rig}}_{v}$
from Theorem \tup{\ref{thm:242}}.
\end{thmstar}
 
\begin{corstar}[\cite{Ye8}] \label{cor:125}
The assignment $A' \mapsto \KK_{A'}$ extends uniquely to a quasi-coherent 
sheaf on the small \'etale site of $\opn{Spec}(A)$. 
\end{corstar}
 
This corollary will be crucial when we get to DM stacks. 

\begin{thmstar}[\cite{Ye8}] \label{thm:130} 
Suppose $u : A \to B$ is a finite homomorphism of EFT $\K$-rings. 
There is a unique nondegenerate backward homomorphism 
$\opn{tr}^{\mrm{res}}_{u} : \mcal{K}_B \to \mcal{K}_{A}$
in $\cat{C}_{\mrm{str}}(A)$, called the {\em rigid trace homomorphism}, such 
that 
$\opn{Q}(\opn{tr}^{\mrm{res}}_{u}) : 
(\KK_B, \rho_B) \to (\KK_{A}, \rho_{A})$
is the rigid trace morphism $\opn{tr}^{\mrm{rig}}_{u}$
from Theorem \tup{\ref{thm:241}}.
\end{thmstar}
 
The uniqueness implies that the rigid trace is functorial on finite 
homomorphisms. 

\begin{thmstar}[\cite{Ye8}] \label{thm:131}
Suppose we are given a finite ring homomorphism $u : A \to  B$ and
an essentially \'etale ring homomorphism $v : A \to A'$.
Define the ring $B' := A' \ot_A B$, and the homomorphisms 
$u' : A' \to B'$ and $w : B \to  B'$ as in Theorem \tup{\ref{thm:240}}, which 
are finite and essentially \'etale, respectively. Then the diagram 
\[ \UseTips \xymatrix @C=6ex @R=5ex {
\KK_A
\ar[d]_{\opn{q}^{\mrm{res}}_{v}}
&
\KK_B
\ar[l]_{\opn{tr}^{\mrm{res}}_u}
\ar[d]^{\opn{q}^{\mrm{res}}_{w}}
\\
\KK_{A'}
&
\KK_{B'}
\ar[l]_{\opn{tr}^{\mrm{res}}_{u'}}
} \]
in $\cat{C}_{\mrm{str}}(A)$ is commutative. 
\end{thmstar}

\begin{thmstar}[\cite{Ye8}] \label{thm:132}
For every homomorphism $u : A \to  B$ in $\cat{Rng} \eftover \K$
there is a unique backward homomorphism of graded $A$-modules 
$\opn{tr}^{\mrm{ires}}_{u} : \mcal{K}_B \to \mcal{K}_{A}$,
called the {\em ind-rigid trace homomorphism}, satisfying these conditions: 
\begin{itemize} 
\item[$\triangleright$] If $u$ is a finite homomorphism, then 
$\opn{tr}^{\mrm{ires}}_{u} = \opn{tr}^{\mrm{res}}_{u}$,
the rigid trace from Theorem \tup{\ref{thm:130}}. 

\item[$\triangleright$] The ind-rigid trace $\opn{tr}^{\mrm{ires}}_{u}$ 
commutes with the rigid \'etale-localizations $\opn{q}^{\mrm{res}}_{v}$ from 
Theorem \tup{\ref{thm:126}}.
 
\item[$\triangleright$] Functoriality: 
$\opn{tr}^{\mrm{ires}}_{u_2 \circ u_1} = \opn{tr}^{\mrm{ires}}_{u_1} \circ \,
\opn{tr}^{\mrm{ires}}_{u_2}$. 
\end{itemize} 
\end{thmstar}

\begin{rem} \label{rem:299}
If $u : A \to B$ is not finite, then usually the ind-rigid trace 
$\opn{tr}^{\mrm{ires}}_{u}$ is not a morphism in $\cat{C}_{\mrm{str}}(A)$, 
i.e.\ it does not commute with the differentials. This same behavior already 
occurred for the trace in \cite[Section VI.4]{RD}. See the Rigid Residue 
Theorem \ref{thm:1} for the case of a proper map of schemes. 
\end{rem}
 
Observe that the rigid auto-duality functor $\opn{D}^{\mrm{rig}}_A$ from 
(\ref{eqn:130}) now becomes 
\[ \opn{D}^{\mrm{rig}}_A = \opn{Hom}_{A}(-, \KK_A) \ : \
\cat{D}^{}_{\mrm{f}}(A)^{\mrm{op}} \to \cat{D}^{}_{\mrm{f}}(A) .  \]
This means that the twisted induction pseudofunctor $\opn{TwInd}$ gets a lot 
more structure from the results of this section. 

A ring homomorphism $v : A \to A'$ is called {\em faithfully \'etale} 
if it is \'etale and faithfully flat. In other words, if the map of affine 
schemes $\opn{Spec}(A') \to \opn{Spec}(A)$ is \'etale and surjective.

\begin{thmstar}[\'Etale Codescent, \cite{Ye8}] \label{thm:133}
Suppose $v : A \to A'$ is a faithfully \'etale ring homomorphism.
Let $A'' := A' \ot_A A'$, and let $w_1, w_2 : A' \to A''$ be the two inclusions.
Then for every integer $i$ the sequence of $A$-module homomorphisms
\[ \KK^i_{A''} 
\xar{ \opn{tr}^{\mrm{ires}}_{w_{\scalesubsub{2}}} - 
\opn{tr}^{\mrm{ires}}_{w_{\scalesubsub{1}}} \, } 
\KK^i_{A'} \xar{ \opn{tr}^{\mrm{ires}}_{v} \, } \KK^i_A \to 0 \]
is exact. 
\end{thmstar}
 
We have not encountered any similar result in the literature. For us this 
theorem will be required for producing ind-rigid traces for maps of DM 
stacks.

\begin{rem} \label{rem:295}
A few words on notation. In Section \ref{sec:RRD-rings} we used the
notation $\opn{tr}^{\mrm{rig}}_{u}$ and $\opn{q}^{\mrm{rig}}_{v}$
for the rigid trace and rigid \'etale-localization morphisms, respectively,
which are morphisms between dualizing complexes in $\cat{D}(A)$.
We use the modified notation $\opn{tr}^{\mrm{res}}_{u}$ and 
$\opn{q}^{\mrm{res}}_{v}$ for the corresponding homomorphisms between rigid 
residue complexes in $\cat{C}_{\mrm{str}}(A)$. 
They are related as follows: 
$\opn{tr}^{\mrm{rig}}_{u} = \opn{Q}(\opn{tr}^{\mrm{res}}_{u})$
and
$\opn{q}^{\mrm{rig}}_{v} = \opn{Q}(\opn{q}^{\mrm{res}}_{v})$,
where $\opn{Q} : \cat{C}_{\mrm{str}}(A) \to \cat{D}(A)$ 
is the categorical localization functor. 
\end{rem}

\begin{rem} \label{rem:260}
It is important to mention that the indecomposable injective $A$-module $J(\p)$ 
is unique only up to isomorphism, and without added structure such an 
isomorphism cannot be made unique or "canonical". To see this, note that 
the automorphism group of $J(\p)$
is the group of invertible elements 
${\what{A}_{\p}}^{\, \times}$ of the complete local ring $\what{A}_{\p}$, and 
the  subgroup of $\what{A}_{\p}^{\, \times}$ consisting of the elements 
congruent to $1$ modulo $\p$ (usually a huge subgroup) acts trivially on the 
submodule $\kk(p) \sub J(\p)$. This implies that there is nothing canonical 
about the embedding $\kk(\p) \inj J(\p)$.

There are ways to enhance $J(\p)$ with some extra structure. 
In our context, the rigid residue complex 
$(\mcal{K}_{A_{\p}}, \rho_{A_{\p}})$ of $A_{\p}$ (relative to $\K$) has the 
property that the $A_{\p}$-module $\KK_{A_{\p}}^{-d}$, where 
$d := \opn{rig{.}dim}_{\K}(\p)$, is isomorphic to $J(\p)$, by Definition 
\ref{dfn:250}. A more global way to "rigidify" $J(\p)$ is to consider the 
$A$-module $\Ga_{\p}(\KK_{A}^{-d})$, the $\p$-torsion submodule of 
$\KK_{A}^{-d}$. The forward functoriality 
$\opn{q}^{\mrm{res}}_{v} : \KK_{A} \to \KK_{A_{\p}}$,
which is associated to the localization ring homomorphism $v : A \to A_{\p}$,
induces an isomorphism 
$\Ga_{\p}(\KK_{A}^{-d}) \iso \KK_{A_{\p}}^{-d}$. 
According to Theorems \ref{thm:126}, \ref{thm:130} and \ref{thm:131}, the 
modules $\KK_{A_{\p}}^{-d}$ and $\Ga_{\p}(\KK_{A}^{-d})$
are functorial, in two directions, in the rings $A_{\p}$ and $A$ respectively. 
Thus these are "canonical incarnations" of $J(\p)$. 

Other ways to obtain "canonical incarnations" of $J(\p)$ are using {\em local 
duality} -- this is implicit throughout \cite{RD}; or using coefficient fields 
and higher completions (see e.g.\ \cite{Ye1} and \cite{Ye2}). 
\end{rem}

%\cleardoublepage
\section{Rigidity, Residues and Duality for Schemes}
 
Recall that $\K$ is a regular finite dimensional noetherian base ring.
We use the notation $\cat{Sch} \ftover \K$ for the category of finite type (FT) 
$\K$-schemes. 

Consider a FT $\K$-scheme $X$. If $U \sub X$ is an affine open set, then 
$A := \Gamma(U, \mcal{O}_X)$ is a FT $\K$-ring.
Given a quasi-coherent $\mcal{O}_X$-module, for every affine open set $U$
we get an $A$-module $M := \Gamma(U, \mcal{M})$.
The resulting functor
$\Gamma(U, -) : \cat{QCoh}(X) \to \cat{Mod}(A)$
is exact. 

\begin{rem} \label{rem:261}
The theory discussed in this section extends easily to {\em essentially finite 
type $\K$-schemes}. But the geometry of EFT $\K$-schemes is a bit messy 
(or perhaps the foundations have not been studied well enough yet), and 
therefore we will not talk about them here. The paper \cite{Ye9} does deal 
with EFT $\K$-schemes, by making special adjustments to overcome a few 
annoying technicalities. 
\end{rem}

If $V \sub U$ is an inclusion of affine open sets in $X$, then geometric 
restriction
$\opn{rest}_{V / U} : \Gamma(U, \mcal{O}_X) \to \Gamma(V, \mcal{O}_X)$
is an {\em \'etale ring homomorphism}. For 
$\MM \in \cat{QCoh}(X)$ we get a nondegenerate forward homomorphism 
$\opn{rest}_{V / U} : \Gamma(U, \mcal{M}) \to \Gamma(V, \mcal{M})$
of $\Gamma(U, \mcal{O}_X)$-modules.

According to the Grothendieck-Matlis classification of injective 
$\OO_X$-modules from \cite[Section II.7]{RD}, every quasi-coherent 
$\OO_X$-module $\MM$ admits a monomorphism $\MM \inj \II$, where $\II$ is an 
$\OO_X$-module that is both quasi-coherent and injective in the category \lb
$\cat{Mod}(X)$ of all $\OO_X$-modules. In particular, the object $\II$ is 
injective in the full subcategory $\cat{QCoh}(X)$ of $\cat{Mod}(X)$. We call 
such an object $\II$ a {\em quasi-coherent injective $\OO_X$-module}. Every 
quasi-coherent injective $\OO_X$-module $\II$ is a direct sum (often infinite) 
of indecomposable quasi-coherent injective $\OO_X$-modules. These 
indecomposable quasi-coherent injective $\OO_X$-modules are classified by 
points $x \in X$. For every point $x$ the corresponding indecomposable 
injective $\JJ(x)$ is a skyscraper quasi-coherent $\OO_X$-module with support 
the closed set $\ol{\{ x \}}$, 
and whose stalk at $x$ is the injective hull of 
the residue field $\kk(x)$ over the local ring $\OO_{X, x}$. To match this with 
the ring theoretic Matlis classification from Section 3, when 
$X = \opn{Spec}(A)$ is affine, and $x = \p$, then $\JJ(x) \cong J(\p)$
as modules over $\OO_{X, x} \cong A_{\p}$.  

\begin{dfn}[\cite{Ye9}] \label{dfn:255}
Let $X$ be a FT $\K$-scheme. A {\em rigid residue complex on $X$} (relative to 
$\K$) is a pair $(\mcal{K}_X, \bsym{\rho}_X)$ consisting of these data:
\begin{itemize}
\rmitem{a} A bounded complex $\mcal{K}_X$ of quasi-coherent injective  
$\mcal{O}_X$-modules,

\rmitem{b} A collection $\bsym{\rho}_X := \{ \rho_U \}$ 
indexed by the affine open sets $U \sub X$, where each 
$\rho_U$ is a rigidifying isomorphism for the complex of  
$\Gamma(U, \mcal{O}_X)$-modules $\Gamma(U, \mcal{K}_X)$, 
in the sense of Definition \ref{dfn:230}.
\end{itemize}
There are two conditions:
\begin{enumerate}
\rmitem{i} For every $U$, the pair 
$\bigl( \Gamma(U, \mcal{K}_X), \rho_U \bigr)$
is a rigid residue complex over the 
ring  $\Gamma(U, \mcal{O}_X)$, in the sense of Definition \ref{dfn:250}.

\rmitem{ii} For every inclusion $V \sub U$ of affine open sets, the 
nondegenerate forward homomorphism 
$\opn{rest}_{V / U} : \Gamma(U, \mcal{K}_X) \to \Gamma(V, \mcal{K}_X)$
is the unique rigid \'etale-localization homomorphism between these rigid 
residue complexes (see Theorem$^{(*)}$ \ref{thm:126}). 
\end{enumerate}
The collection $\bsym{\rho}_X$ is called a {\em rigid structure} on the complex 
$\KK_X$. 
\end{dfn}

\begin{dfn}[\cite{Ye9}] \label{dfn:256}
Suppose $(\mcal{K}_X, \bsym{\rho}_X)$ and $(\mcal{K}'_X, \bsym{\rho}'_X)$
are two rigid residue complexes on the FT $\K$-scheme $X$. 
A {\em morphism of rigid residue complexes}
$\phi : (\mcal{K}_X, \bsym{\rho}_X) \to (\mcal{K}'_X, \bsym{\rho}'_X)$
is a homomorphism $\phi : \KK_X \to \KK'_X$ of complexes of
$\mcal{O}_X$-modules, such that for every
affine open set $U$, with $A := \Gamma(U, \mcal{O}_X)$, the 
homomorphism $\Gamma(U, \phi)$ is a morphism in 
$\cat{C}(A)_{\mrm{res} / \K}$.
\end{dfn}

The next four theorems are easy consequences of the corresponding theorems for 
rings in Section \ref{sec:rrr-rings}, combined with descent for QC sheaves on 
schemes. 
 
\begin{thmstar}[\cite{Ye9}] \label{thm:135}
Every FT $\K$-scheme $X$ has a rigid residue complex 
$(\mcal{K}_X, \bsym{\rho}_X)$, and it is unique up to a unique isomorphism of 
rigid residue complexes. 
\end{thmstar}

\begin{thmstar}[\cite{Ye9}] \label{thm:136}
For every map $f : Y \to X$ between FT $\K$-schemes there is a unique 
homomorphism of {\em graded quasi-coherent $\mcal{O}_X$-modules}
$\opn{tr}^{\mrm{ires}}_f : f_*(\mcal{K}_Y) \to \mcal{K}_X$,
called the {\em ind-rigid trace homomorphism}, 
which extends the ind-rigid trace homomorphism $\opn{tr}^{\mrm{ires}}_{(-)}$
on $\K$-rings from Theorem \tup{\ref{thm:132}}. The homomorphism 
$\opn{tr}^{\mrm{ires}}_f$ is functorial in $f$.   
\end{thmstar}
 
If $f$ is a {\em finite} map of schemes, then $\opn{tr}^{\mrm{ires}}_f$ 
is a homomorphism of complexes -- this is immediate from the result for rings. 
The case of a {\em proper map} will be mentioned soon. 

\begin{thmstar}[\cite{Ye9}] \label{thm:137}
For every \'etale map $g : X' \to X$ between FT $\K$-schemes
there is a unique homomorphism of complexes of quasi-coherent 
$\mcal{O}_X$-modules
$\opn{q}^{\mrm{res}}_g : \mcal{K}_X \to g_*(\mcal{K}_{X'})$,
called the {\em rigid \'etale-localization homomorphism}, 
which extends the rigid \'etale-localiza\-tions homomorphism 
$\opn{q}^{\mrm{res}}_{(-)}$ on $\K$-rings from Theorem \tup{\ref{thm:126}}. 
The homomorphism $\opn{q}^{\mrm{res}}_g$ is functorial in such maps $g$.  
\end{thmstar}
 
\begin{thmstar}[\cite{Ye9}] \label{thm:138}
The ind-rigid trace homomorphism $\opn{tr}^{\mrm{ires}}_{(-)}$ and the 
rigid \'etale-locali\-zation homomorphism $\opn{q}^{\mrm{res}}_{(-)}$ 
commute with each other, as in Theorem$^{(*)}$ \tup{\ref{thm:131}} but 
geometrized.
\end{thmstar}

The next result is really just a geometric version of Theorem \ref{thm:121}, 
without any significant difficulties. For a scheme $X$ with rigid residue 
complex $(\KK_A, \rho_A)$, we define
\begin{equation} \label{eqn:295}
\opn{D}_{X}^{\mrm{rig}} := \Hom_{X}(-, \KK_X) : \cat{D}(X) \to \cat{D}(X)  .
\end{equation}
Here $\cat{D}(X)$ is the derived category of $\OO_X$-modules.
As in the affine case, there are equivalences 
$\opn{D}_{X}^{\mrm{rig}} : \cat{D}^{\star}_{\mrm{c}}(X) \to 
\cat{D}^{-\star}_{\mrm{c}}(X)$
for reversed boundedness conditions $\star$ and $-\star$. 

\begin{thm}[{\cite{YZ2}, \cite{Ye9}}] \label{thm:290} 
There is a unique pseudofunctor 
\[ \opn{TwInvIm} : (\cat{Sch} \ftover \K)^{\mrm{op}} \to \cat{TrCat} \over \K \]
called {\em twisted inverse image}, with these properties\hsp{0.1}\tup{:}
\begin{enumerate} 
\item To an object $X \in \cat{Sch} \ftover \K$ it assigns the triangulated 
category $\cat{D}^{+}_{\mrm{c}}(X)$.
 
\item To a morphism $f : Y \to X$ in $\cat{Sch} \ftover \K$ it assigns the 
triangulated functor 
\[ \opn{TwInvIm}_{f} = f^! \ := \ 
\opn{D}_{Y}^{\mrm{rig}} \circ \, \mrm{L} f^* \circ \opn{D}_{X}^{\mrm{rig}}
\ :  \ \cat{D}^{+}_{\mrm{c}}(X) \to \cat{D}^{+}_{\mrm{c}}(Y) . \]
 
\item The composition isomorphisms of \, $\opn{TwInvIm}$ come from the 
evaluation isomorphisms of the rigid auto-duality functors 
\tup{(\ref{eqn:295})}, combined 
with the composition isomorphisms of the pseudofunctor $f \mapsto \mrm{L} f^*$.
\end{enumerate}
\end{thm}

The functorial properties of the rigid residue complexes (stated in the 
theorems above and below) provide the twisted inverse image pseudofunctor with 
more structure. 

Here is the rigid version of \cite[Theorem VII.2.1]{RD}.

\begin{thmstar}\tup{(Rigid Residue Theorem, \cite{Ye9})} \label{thm:1}
Let $f : Y \to X$ be a proper map between FT $\K$-schemes. Then the ind-rigid 
trace
$\opn{tr}^{\mrm{ires}}_f  : f_*(\mcal{K}_Y) \to \mcal{K}_X$
is a homomorphism of complexes. 
\end{thmstar}

The idea of the proof (imitating \cite{RD}) is to reduce to the case when 
$X = \opn{Spec}(A)$ for a local artinian ring 
$A \in \cat{Rng} \eftover \K$, and $Y = \mbf{P}^1_A$, the projective line.
We do it using Theorems \ref{thm:136}, \ref{thm:137} and \ref{thm:138}.
Then we do a calculation of residues, which similar to the one in the proof of 
the  classical Residue Theorem for $\mbf{P}^1(\mbb{C})_{\mrm{an}}$. The only 
minor complication is that a rational differential form on $Y$ relative to $A$, 
i.e.\ a form $\om \in \Om^1_{A(t) / A}$, where $A(t)$ is the total ring of 
fractions of $Y$, might have poles at closed points of $Y$ that do not belong 
to $Y(A)$. But the poles of $\om$ always belong to $Y(A')$ for a suitable 
finite faithfully flat $A$-ring $A'$, and by performing the base change 
$A \to A'$ the calculation can proceed. 

The next theorem is the rigid version of \cite[Theorem VII.3.3]{RD}.

\begin{thmstar}\tup{(Rigid Duality Theorem, \cite{Ye9})} \label{thm:2}
Let $f : Y \to X$ be a proper map between FT $\K$-schemes. 
For every $\mcal{N} \in \cat{D}^{}_{\mrm{qc}}(Y)$ the morphism  
\[ \DR f_* \bigl( \DR \Hom_{Y}(\mcal{N}, \mcal{K}_Y) \bigr) \to
\DR \Hom_{X} \bigl( \DR f_* (\mcal{N}), \mcal{K}_X \bigr) \]
in $\cat{D}(X)$, which is induced by the ind-rigid trace
$\opn{tr}^{\mrm{ires}}_f : f_*(\mcal{K}_Y) \to \mcal{K}_X$,
is an isomorphism. 
\end{thmstar}

The proof of  imitates the proof of the corresponding
theorem in \cite{RD}, once we have the Theorem$^{(*)}$ \ref{thm:1} at hand.

Theorem$^{(*)}$ \ref{thm:1} implies that for every complex 
$\MM \in \cat{D}^{+}_{\mrm{c}}(X)$ there is an ind-rigid trace morphism 
\begin{equation} \label{eqn:285}
\opn{tr}^{\mrm{irig}}_{f, \MM} : 
\mrm{R} f_*(f^!(\MM)) \to \MM  
\end{equation}
in $\cat{D}(X)$, which is functorial in $\MM$. 
Here is our version of \cite[Corollary 3.4]{RD}. 

\begin{corstar}\tup{(\cite{Ye9})} \label{cor:285}
Let $f : Y \to X$ be a proper map between FT $\K$-schemes. 
For every $\mcal{M} \in \cat{D}^{+}_{\mrm{c}}(X)$ and 
$\mcal{N} \in \cat{D}^{-}_{\mrm{qc}}(Y)$ the morphism  
\[ \DR f_* \bigl( \DR \Hom_{Y} \bigl( \mcal{N}, f^!(\MM) \bigr) \bigr) \to
\DR \Hom_{X} \bigl( \DR f_* (\mcal{N}), \mcal{M} \bigr) \]
in $\cat{D}(X)$, which is induced by the ind-rigid trace
$\opn{tr}^{\mrm{irig}}_{f, \MM}$, is an isomorphism. 
\end{corstar}

\begin{rem}
At this stage we have recovered almost all the results of the original book 
\cite{RD}. One advantage of our rigidity approach is that it is much cleaner
and shorter than the original approach in \cite{RD}. We can also treat EFT maps 
of schemes. Another advantage of the rigidity approach, as we shall see next, 
is that  it gives rise to a useful duality theory for DM stacks. 
On the down side, we must work with a fixed regular base ring $\K$. 

For complements to the book \cite{RD}, and for alternative approaches, due to 
J. Lipman, B. Conrad, A. Neeman, the author of this article, and others, please 
consult \cite{Co}, \cite{LH}, \cite{Ne1}, \cite{Ne2}, \cite{Ye7} and 
their references. 
\end{rem}

%\cleardoublepage
\section{Rigidity, Residues and Duality for DM Stacks}

Here is a brief recollection of facts about Deligne-Mumford (DM) stacks and 
quasi-coherent sheaves on them, extracted from \cite{SP} and \cite{Ol}, with 
some additional  elaboration. This will help us state our results. 

As a first approximation, it useful to think of a DM stack $\XX$ as a 
scheme equipped with an extra structure: the points of $\XX$ are bunched 
together as objects of a groupoid, in an intricate geometric way. 
(This can be understood, for instance, as the way $\XX$ sits above its coarse 
moduli space $X$; we will return to this idea later.)

To give a more precise description of DM stacks,
we need to place ourselves in the realm of sheaves and stacks 
on $(\cat{Sch} \ftover \K)_{\mrm{et}}$, the big \'etale site of FT $\K$-schemes.

For us a {\em prestack of groupoids} on $\cat{Sch} \ftover \K$ is 
a pseudofunctor 
$\XX : (\cat{Sch} \ftover \K)^{\mrm{op}} \to \cat{Grpd}$,
where $\cat{Grpd}$ is the $2$-category of groupoids.
The prestack $\XX$ is a {\em stack} if it satisfies descent for morphisms and 
descent for objects with respect to coverings in 
$(\cat{Sch} \ftover \K)_{\mrm{et}}$. 
See \cite[Section 2]{Ye4} for details on stacks, as viewed from the 
pseudofunctor approach. 

Most algebraic geometry texts (including \cite{Ol} and \cite{SP}) prefer 
viewing a prestack of groupoids on $\cat{Sch} \ftover \K$ as a {\em category 
fibered in groupoids}. The passage from this approach to ours is this: 
a prestack of groupoids $\XX$ in our sense is the same as a category 
fibered in groupoids equipped with a cleavage. See \cite[Remark 2.5]{Ye4}.

A benefit of our approach to stacks of groupoids is that it lets us 
refer to a sheaf of sets as a very trivial kind of stack. Indeed, a sheaf of 
sets $X$ on $(\cat{Sch} \ftover \K)_{\mrm{et}}$ can be seen as a stack of 
groupoids, in which the only local isomorphisms between local objects are the 
identity automorphisms. 

We know how a scheme $X$ becomes a sheaf of sets on 
$(\cat{Sch} \ftover \K)_{\mrm{et}}$ -- it is the Yoneda incarnation of $X$ 
as the sheaf $\Hom(-, X)$. An {\em algebraic space} $X$ is a more 
complicated sheaf of sets on $(\cat{Sch} \ftover \K)_{\mrm{et}}$; 
it is not representable, but only "locally representable". One way to express 
this condition is that the sheaf $X$ is 
isomorphic to the quotient sheaf $U_0 / U_1$ of an {\em \'etale equivalence 
relation in schemes} $(U_1 \rightrightarrows U_0)$. Let's refer to the 
isomorphism of sheaves $p : U_0 / U_1 \iso X$ as a {\em presentation} of 
$X$. See \cite[Proposition 5.2.5]{Ol}.

There is an analogous way to say when an abstract stack $\XX$
on $(\cat{Sch} \ftover \K)_{\mrm{et}}$ is a DM stack. 
The condition is that $\XX$ is equivalent, in the $2$-category of stacks of 
groupoids on $(\cat{Sch} \ftover \K)_{\mrm{et}}$,
to the quotient stack $[X_0 / X_1]$ of an {\em \'etale groupoid in algebraic 
spaces} \lb $(X_1 \rightrightarrows X_0)$. Note that the groupoid 
$(X_1 \rightrightarrows X_0)$ involves a third algebraic space $X_2$, with a map
$X_2 \to X_1 \times_{X_0} X_1$ that encodes composition of morphisms. 
Again we call the equivalence
$p : [X_0 / X_1] \equiv \XX$ a presentation of $\XX$. See 
\cite[Theorem \href{https://stacks.math.columbia.edu/tag/04TK}{\tt tag=04TK}
and Lemma \href{https://stacks.math.columbia.edu/tag/04T5}{\tt tag=04T5}]{SP}.
We shall use the term {\em map of stacks} for a $1$-morphism
$f : \YY \to  \XX$ between stacks of groupoids on 
$(\cat{Sch} \ftover \K)_{\mrm{et}}$. 

Thus there are four stages, or levels of complexity, in the geometry of 
FT DM $\K$-stacks: 
\begin{equation} \label{280}
\cat{AfSch} \ftover \K 
\sub \cat{Sch} \ftover \K 
\sub \cat{AlgSp} \ftover \K 
\sub \cat{DMStk} \ftover \K  .
\end{equation}
These are the affine schemes, schemes, algebraic spaces and DM stacks, all of 
FT over $\K$. If we ignore the fact that $\cat{DMStk} \ftover \K$ is 
a $2$-category (i.e.\ if we 
forget the $2$-morphisms between $1$-morphisms of DM stacks), formula 
(\ref{280}) describes full embeddings of categories. 
Morphisms in these categories will be called {\em maps}. 
In what follows all our stacks (and algebraic spaces and schemes) 
are FT over the base ring $\K$; and the rings are EFT over $\K$. 

A map $f : \YY \to \XX$ of DM stacks is {\em surjective} if it is locally 
essentially surjective on objects. (This is a property of a map between 
abstract stacks of groupoids, not just of algebraic stacks.) 

A map $f : \YY \to \XX$ between DM stacks is called {\em \'etale} if there exist
presentations 
$p : [X_0 / X_1] \equiv \XX$ and $q : [Y_0 / Y_1] \equiv \YY$,
and an \'etale map of algebraic spaces
$g : Y_0 \to X_0$, such that the diagram 
\begin{equation} \label{eqn:270}
\UseTips \xymatrix @C=7ex @R=5ex {
Y_0
\ar[r]^{q}
\ar[d]_{g}
&
\YY
\ar[d]^{f}
\\
X_0
\ar[r]^{p}
&
\XX }
\end{equation}
of stacks is commutative up to isomorphism.
For this condition we need to know what is an \'etale map 
$g : Y_0 \to X_0$ between algebraic spaces, and this is done by going down the 
staircase of embeddings (\ref{280}).
(This same yoga is how other local properties of a map of DM 
stacks $f : \YY \to \XX$, such as being smooth or flat, are defined -- they 
are reflected by the corresponding properties of $g : Y_0 \to X_0$.)
The more conventional criterion for $f : \YY \to \XX$ to be 
\'etale is in \cite[Definition 8.2.6]{Ol}.

Presentations of algebraic spaces and DM stacks allow us to describe 
{\em quasi-coherent sheaves} on these geometric objects. Thus, a presentation 
$p : U_0 / U_1 \iso X$ of an algebraic space $X$ induces an equivalence of 
abelian categories 
\begin{equation} \label{eqn:211}
p^* : \cat{QCoh}(X) \to \cat{QCoh} (U_1 \rightrightarrows U_0) , 
\end{equation}
where the latter refers to the category of equivariant quasi-coherent 
sheaves on the equivalence relation $(U_1 \rightrightarrows U_0)$. 
Analogously, a presentation $p : [X_0 / X_1] \equiv \XX$ of a DM stack 
$\XX$ induces an equivalence of abelian categories 
\begin{equation} \label{eqn:212}
p^* : \cat{QCoh}(\XX) \to \cat{QCoh} (X_1 \rightrightarrows X_0) , 
\end{equation}
where now $(X_1 \rightrightarrows X_0)$ is an \'etale groupoid in algebraic 
spaces. See \cite[Sections 7.1 and 9.2]{Ol} and 
\cite[Proposition \href{https://stacks.math.columbia.edu/tag/03M3}{\tt 
tag=03M3} and Lemma 
\href{https://stacks.math.columbia.edu/tag/06WV}{\tt tag=06WV}]{SP}.

The discussion above suggests a mechanism for producing a quasi-coherent 
sheaf $\MM_{\XX}$ on a FT DM $\K$-stack $\XX$. Consider a presentation 
$[X_0 / X_1] \equiv \XX$, the finitely many schemes $U$ involved in presenting 
the algebraic spaces $X_i$, the finitely many affine schemes $\opn{Spec}(A)$
involved in suitable hypercoverings of these schemes, and the finitely many 
\'etale maps between these geometric objects,  going down the staircase in 
the hierarchy (\ref{280}). 
The first step in this mechanism is to provide an $A$-module $M_A$ 
for every affine scheme $\opn{Spec}(A)$ involved, and a consistent collection
of nondegenerate forward homomorphisms $M_A \to M_{A'}$ for all the \'etale ring
homomorphisms $A \to A'$ involved. The second step is to glue these modules
into QC sheaves $\MM_U$ on all the schemes $U$ involved. There will be 
an induced consistent collection of isomorphisms $g^*(\MM_U) \iso \MM_{U'}$ for 
the \'etale maps $g : U' \to  U$ involved. The third step is to use 
the equivalences (\ref{eqn:211}) associated to presentations to obtain 
QC sheaves $\MM_{X_i}$ on the algebraic spaces $X_i$, with consistent pullback 
isomorphisms. Finally, formula (\ref{eqn:212}) will give us the 
desired quasi-coherent sheaf $\MM_{\XX}$. 

Suppose we have a map $f : \YY \to \XX$ of DM stacks, and QC sheaves
$\MM \in  \cat{QCoh}(\XX)$ and $\NN \in  \cat{QCoh}(\YY)$. 
To specify a homomorphism 
$\phi : f_{*}(\NN) \to \MM$ in $\cat{QCoh}(\XX)$
we can use presentations. After having arranged for 
presentations $p : [X_0 / X_1] \equiv \XX$ and $q : [Y_0 / Y_1] \equiv \YY$,
and a map of algebraic spaces $\til{f}_0 : Y_0 \to X_0$, such that the diagram 
(\ref{eqn:270}), with $\til{f}_0$ instead of $g$,  is commutative (up to 
isomorphism), the map $\til{f}_0$ can be extended to a map 
\begin{equation} \label{eqn:284}
\til{f} = (\til{f}_0, \til{f}_1, \til{f}_2) : 
(Y_1 \rightrightarrows Y_0) \to (X_1 \rightrightarrows X_0)
\end{equation}
of groupoids in algebraic spaces (see
\cite[
Definition \href{https://stacks.math.columbia.edu/tag/043W}{\tt tag=043W}
and 
Section  \href{https://stacks.math.columbia.edu/tag/04TJ}{\tt tag=04TJ}]{SP}).
By pullback along $p$ and $q$ we get quasi-coherent sheaves $\MM_i \in 
\cat{QCoh}(X_i)$ and $\NN_i \in \cat{QCoh}(Y_i)$ satisfying the appropriate 
relations. Then we need to arrange for homomorphisms 
$\til{\phi}_i : \til{f}_{i \, *}(\NN_i) \to \MM_i$
that respect the relations; and these will determine $\phi$ by the equivalence 
(\ref{eqn:212}). As before, this construction requires going 
down the staircase of the hierarchy (\ref{280}). 

After recalling and elucidating the mechanism of gluing quasi-coherent sheaves 
on DM stacks, we can finally talk about rigid residue complexes. As we shall 
see, the definition of a rigid residue complex on a DM stack $\XX$ is
very similar to the definition for a scheme, with only a minor adjustment. (We 
will skip the algebraic space case, since it is identical to the DM 
stack case.)

\begin{dfn}[\cite{Ye10}] \label{dfn:270}
Let $\XX$ be a FT DM $\K$-stack. A {\em rigid residue complex on $\XX$} 
(relative to $\K$) is a pair $(\mcal{K}_{\XX}, \bsym{\rho}_{\XX})$ consisting 
of these data:
\begin{itemize}
\rmitem{a} A bounded complex $\mcal{K}_X$ of quasi-coherent injective  
$\mcal{O}_{\XX}$-modules.

\rmitem{b} A collection 
$\bsym{\rho}_{\XX} := \{ \rho_{(U, g)} \}$ 
indexed by \'etale maps $g : U \to \XX$ from affine $\K$-schemes $U$,
where each $\rho_{(U, g)}$ is a rigidifying isomorphism for the complex of  
$\Gamma(U, \mcal{O}_U)$-modules $\Gamma(U, g^*(\mcal{K}_{\XX}))$, 
in the sense of Definition \ref{dfn:230}.
\end{itemize}
There are two conditions:
\begin{enumerate}
\rmitem{i} For every \'etale map $g : U \to  \XX$ from an affine scheme 
$U = \opn{Spec}(A)$, letting
$\KK_A := \Gamma(U, g^*(\mcal{K}_{\XX}))$, 
the pair $(\KK_A, \rho_{(U, g)})$ 
is a rigid residue complex over the ring $A$, in the sense of Definition 
\ref{dfn:250}.

\rmitem{ii} For every commutative diagram 
\[ \UseTips \xymatrix @C=7ex @R=4ex {
U_2
\ar[r]^{h}
\ar[dr]_{g_2}
&
U_1
\ar[d]^{g_1}
\\
&
\XX
} \]
of \'etale maps, where  $U_1$ and $U_2$ are affine schemes, 
the homomorphism $h^*$ arising from the equality $g_2 = g_1 \circ h$
is the unique nondegenerate rigid \'etale-localization homomorphism
\[ h^* : \bigl( \Gamma(U_1, g_1^*(\mcal{K}_{\XX})), \rho_{(U_1, g_1)} \bigr) \to
\bigl( \Gamma(U_2, g_2^*(\mcal{K}_{\XX})), \rho_{(U_2, g_2)} \bigr) , \]
in the sense of Theorem$^{(*)}$ \ref{thm:126}.
\end{enumerate}
The collection $\bsym{\rho}_{\XX}$ is called a {\em rigid structure} on the 
complex $\KK_{\XX}$. 
\end{dfn}
 
\begin{thmstar} \tup{ (\cite{Ye10})} \label{thm:3}
Let $\XX$ be a FT DM $\K$-stack. There exists a rigid residue complex 
$(\KK_{\XX}, \bsym{\rho}_{\XX})$, and it is unique up to a unique rigid 
isomorphism.
\end{thmstar}

The proof is by applying the mechanism of gluing quasi-coherent sheaves that 
was explained above, where in the first 
step of affine schemes we use Corollary$^{(*)}$ \ref{cor:125}. 
 
\begin{thmstar}\tup{ (\cite{Ye10})} \label{thm:4}
Let $f : \YY \to \XX$ be a map between FT DM $\K$-stacks. 
There is a homomorphism of {\em graded quasi-coherent 
$\OO_{\XX}$-modules} 
$\opn{tr}^{\mrm{ires}}_f : f_*(\mcal{K}_{\YY}) \to \mcal{K}_{\XX}$
called the {\em ind-rigid trace}, extending the ind-rigid trace on FT 
$\K$-rings. 
\end{thmstar}
 
To prove this theorem we apply the construction explained above for 
homomorphisms between quasi-coherent sheaves, starting with 
Theorem$^{(*)}$ \ref{thm:133} for affine schemes. 

\begin{thmstar}[{\cite{Ye10}}] \label{thm:287} 
There is a pseudofunctor 
\[ \opn{TwInvIm} : (\cat{DMStk} \ftover \K)^{\mrm{op}} \to \cat{TrCat} \over \K
, \quad \XX \mapsto \cat{D}^{+}_{\mrm{c}}(\XX) 
, \quad f \mapsto f^!  \]
called the {\em twisted inverse image}, with properties like in Theorem 
\tup{\ref{thm:290}}.
\end{thmstar}

In this theorem we consider $\cat{DMStk} \ftover \K$ as a $1$-category, 
ignoring the $2$-isomorphisms between $1$-morphisms, as we had done in the 
hierarchy (\ref{280}). 
However, the $2$-category structure of $\cat{DMStk} \ftover \K$ is reflected in 
this property: if $f_0, f_1 : \YY \to \XX$ are maps in 
$\cat{DMStk} \ftover \K$
for which a $2$-isomorphism $\ga : f_0 \twoiso f_1$ exists, then the 
triangulated functors 
$f_0^!, f_1^! : \cat{D}^{+}_{\mrm{c}}(\XX) \to \cat{D}^{+}_{\mrm{c}}(\YY)$
are isomorphic. This implies that when $f : \YY \to \XX$ is an 
equivalence of stacks, the functor $f^!$ will also be an equivalence. 

The obvious question now is this: {\em Do the Rigid Residue Theorem and the 
Rigid Duality Theorem hold for a proper map $f : \YY \to \XX$ 
between DM stacks}~{\hspace{-0.2ex}}?
We only know a partial answer. 

By the Keel-Mori Theorem \cite[Theorem 11.1.2]{Ol}, a separated FT $\K$-stack 
$\YY$ has a {\em coarse moduli space} $\pi : \YY \to Y$. The map $\pi$ is 
proper and quasi-finite, and $Y$ is, in general, an {\em algebraic 
space}. 

\begin{dfn} \label{dfn:271}
\begin{enumerate}
\item Let $\YY$ be a separated FT DM $\K$-stack, with coarse moduli space
$\pi : \YY \to Y$. We call $\YY$ a {\em coarsely schematic stack} if its coarse 
moduli space $Y$ is a scheme. 

\item Let $f : \YY \to \XX$ be a separated map between FT DM $\K$-stacks. 
We say that $f$ is a {\em coarsely schematic map} if for some surjective 
\'etale 
map $U \to \XX$ from an affine scheme $U$, the stack 
$\YY' := \YY \times_{\XX} U$ is  coarsely schematic.
\end{enumerate} 
\end{dfn}

This appears to be a rather mild restriction: most DM stacks that come up in
examples are of this kind. See Remark \ref{rem:275} about the necessity of this 
condition.

\begin{thmstar}[Rigid Residue Theorem, \cite{Ye10}] \label{thm:10}
Suppose $f : \YY \to \XX$ is a proper coarsely schematic map of FT DM 
$\K$-stacks. Then the ind-rigid trace
$\opn{tr}^{\mrm{ires}}_f : f_*(\mcal{K}_{\YY}) \to \mcal{K}_{\XX}$
is a homomorphism of complexes of $\OO_{\XX}$-modules. 
\end{thmstar}

It is not expected that duality will hold in this generality. In fact, there
are easy counter examples. The problem is not geometric, but rather a basic 
difficulty with {\em representations of finite groups in positive 
characteristics}.
  
Following \cite{AOV} (see also \cite[Definition 11.3.2]{Ol}), a separated FT DM 
stack $\YY$ is called {\em tame} if for every field 
$L \in \cat{Rng} \eftover \K$, the automorphism groups in the groupoid 
$\YY(L)$, which are always finite, have orders prime to the characteristic of 
$L$. 

\begin{dfn} \label{dfn:275}
A separated map of FT DM $\K$-stacks $f : \YY \to \XX$ is called a {\em tame 
map} if for some surjective \'etale map $U \to \XX$ from an affine scheme $U$, 
the stack $\YY' := \YY \times_{\XX} U$ is tame.
\end{dfn}

\begin{thmstar}[Rigid Duality Theorem, \cite{Ye10}] \label{thm:11}
Suppose $f : \YY \to \XX$ is a proper tame coarsely schematic map between FT
DM $\K$-stacks. 
For every $\mcal{N} \in \cat{D}^{-}_{\mrm{qc}}(\YY)$ the morphism  
\[ \DR f_* \bigl( \DR \Hom_{\YY}(\mcal{N}, \mcal{K}_{\YY}) \bigr) \to
\DR \Hom_{\XX} \bigl( \DR f_* (\mcal{N}), \mcal{K}_{\XX} \bigr) \]
in $\cat{D}(\XX)$, which is induced by the ind-rigid trace
$\opn{tr}^{\mrm{ires}}_f : f_*(\mcal{K}_{\YY}) \to \mcal{K}_{\XX}$,
is an isomorphism. 
\end{thmstar}

Theorem$^{(*)}$ \ref{thm:10} implies that for every complex 
$\MM \in \cat{D}^{+}_{\mrm{c}}(\XX)$ there is an ind-rigid trace morphism 
$\opn{tr}^{\mrm{irig}}_{f, \MM} : 
\mrm{R} f_*(f^!(\MM)) \to \MM$  
in $\cat{D}(\XX)$, which is functorial in $\MM$. 

\begin{corstar}\tup{(\cite{Ye10})} \label{cor:286}
Suppose $f : \YY \to \XX$ is a proper tame coarsely schematic map between FT
DM $\K$-stacks. For every $\mcal{M} \in \cat{D}^{+}_{\mrm{c}}(\XX)$ and 
$\mcal{N} \in \cat{D}^{-}_{\mrm{qc}}(\YY)$ the morphism  
\[ \DR f_* \bigl( \DR \Hom_{\YY}(\mcal{N}, f^!(\MM)) \bigr) \to
\DR \Hom_{\XX} \bigl( \DR f_* (\mcal{N}), \mcal{M} \bigr) \]
in $\cat{D}(\XX)$, which is induced by the ind-rigid trace
$\opn{tr}^{\mrm{irig}}_{f, \MM}$, is an isomorphism. 
\end{corstar}

Here is a sketch of the proofs of Theorems$^{(*)}$ \ref{thm:10} and 
\ref{thm:11}, which we are going to refer to as "Residue" and "Duality", 
respectively. 
Take a surjective \'etale map $U \to \XX$ from an affine scheme 
$U$ such that the stack $\YY' := \YY \times_{\XX} U$ is coarsely schematic.
Consider the commutative up to isomorphisms diagram of maps of stacks
\begin{equation*}
\UseTips \xymatrix @C=5ex @R=3ex {
&
\YY'
\ar[dl]_{\pi'}
\ar[dd]^{f'}
\ar[r]
&
\YY
\ar[dd]^{f}
\\
Y'
\ar[dr]_{g'}
\\
&
U
\ar[r]{}
&
\XX
}
\end{equation*}
where $f'$ is gotten from $f$ by base change, and 
$Y'$ is the coarse moduli space of $\YY'$.
The maps $\pi'$ and $g'$ are both proper.  

Because the functor 
$\cat{QCoh}(\XX) \to \cat{QCoh}(U)$ is faithful, it suffices to prove 
``Residue'' and ``Duality'' for the map $f'$.  
Now $Y'$ is a scheme, so the proper map $g' : Y' \to U$ satisfies both
``Residue'' and ``Duality'', by Theorems$^{(*)}$ \ref{thm:1} and 
\ref{thm:2}.

It remains to verify ``Residue'' and ``Duality'' for the proper map 
$\pi' : \YY' \to Y'$. These properties are \'etale local on $Y'$.
Namely let $V'_1, \ldots, V'_n$ be affine schemes, and let
\begin{equation}  \label{eqn:11}
\coprod\nolimits_{i} V'_i \ \to \ Y'
\end{equation}
be a surjective \'etale map. For any $i$ let 
$\YY'_i := \YY' \times_{Y'} V'_i$.
We get this commutative (up to isomorphism) diagram of maps of stacks
\begin{equation*}
\UseTips \xymatrix @C=6ex @R=5ex {
\coprod_{i} \YY'_i 
\ar[r]
\ar[d]_{\coprod_{i}  \pi'_i}
&
\YY'
\ar[d]^{\pi'}
\\
\coprod_{i} V'_i
\ar[r]
&
Y'
} 
\end{equation*}
It is enough to check ``Residue'' and ``Duality'' for the maps
$\pi'_i : \YY'_i \to V'_i$.
Note that the affine scheme $V'_i$ is the coarse moduli space of the stack 
$\YY'_i$. 

According to \cite[Theorem 11.3.1]{Ol} it is possible to choose a covering 
(\ref{eqn:11}) such that 
\begin{equation}
\YY'_i \approx [W_i / G_i] \quad \tup{and} \quad
V'_i \cong W_i / G_i . 
\end{equation}
Here  $W_i$ is an affine scheme, $G_i$ is a finite group acting on $W_i$,  
$[W_i / G_i]$ is the quotient stack, and $W_i / G_i$ is the quotient scheme. 
Moreover, in the tame case we can assume that 
the order of the group $G_i$
is invertible in the ring $\Ga(W_i, \OO_{W_i})$.

We have now reduced the problem to proving ``Residue'' and ``Duality'' for a
map of stacks 
\[ \pi \, : \, \YY = [W / G] \, \to \, V = W / G \, , \]
where $W = \opn{Spec}(A)$ for some ring $A$, $G$ is a finite group acting on
$A$, and in the tame case the order of $G$ is invertible in $A$. 

Let $A^{G}$ be the ring of $G$-invariants inside $A$.
So $A^{G}$ is a FT $\K$-ring, and the quotient 
scheme is $W / G = \opn{Spec}(A^{G})$. 
Consider the {\em skew group ring} $G \ltimes A$; this is the noncommutative 
(NC) ring appearing in \cite[Exercise 9.E]{Ol}. The NC ring $G \ltimes A$ is 
finite over its center $A^{G}$. It is known that 
$\cat{QCoh}(\YY) \approx \cat{Mod}(G \ltimes A)$
as abelian categories, and the functor $\pi_*$ on the geometric side 
goes to the $G$-invariants functor $(-)^G$ on the NC algebraic side. 
 
We finish the proof by using GD for NC rings, see \cite[Chapters 14-18]{Ye7}. 
 
\begin{rem} \label{rem:275}
It is very likely that the ``coarsely schematic'' condition can be removed 
from Theorems$^{(*)}$ \ref{thm:10} and \ref{thm:11}.
This will probably require a better understanding of the geometry of 
algebraic spaces. 
\end{rem}

\begin{rem} \label{rem:276}
We do not know how to handle {\em Artin stacks}. This is because we do 
not understand the {\em smooth functoriality of the squaring operation}. 
\end{rem}

\begin{rem}
Let us end the article with a brief discussion of related work. 
\begin{enumerate}
\item There is a preprint \cite{Ni} by F. Nironi from 2008 on GD for DM stacks.
It has neither been updated nor published. 
The details in that paper are not clear to us.

\item J. Hall and D. Rydh wrote several papers dealing with derived categories 
on stacks. In their paper \cite{HR} from 2017 they prove a version of GD for a 
finite faithfully flat representable map $f : \YY \to \XX$ of Artin stacks. 

\item A. Neeman wrote the paper \cite{Ne2}, last updated in 2017.
It deals with GD in great generality, and also in a very abstract manner. 
The geometric objects considered in it range from schemes to Artin stacks; and 
the maps between them are required to be concentrated and compactifiable.  
For any such map $f : \YY \to \XX$ the functor $\mrm{R} f_*$ is shown to 
have a right adjoint $f^{\times}$; and under certain further conditions one has 
$f^{\times} \cong f^!$. According to Neeman (in a private communication), his 
work implies that our Corollary \ref{cor:286} holds, when $f : \YY \to \XX$ is 
a concentrated proper map between DM stacks. 
\end{enumerate}

There might be other work on global GD, perhaps also from the direction of 
derived algebraic geometry, but we are not aware of any at present. 
\end{rem}

%\cleardoublepage

\end{document}